\documentclass[12pt]{amsart}

\usepackage{amssymb,amsfonts,amsmath,amsthm,epsfig,tikz,bm,color}

\makeatletter
\def\@seccntformat#1{\csname the#1\endcsname. }
\def\@biblabel#1{#1.}

\makeatother
  
\def\Rm#1{\lowercase\expandafter{\romannumeral#1}}

   \def\Ga{\Gamma}
  \def\D{\Delta} 
  
\def\a{\alpha} \def\b{\beta} \def\g{\gamma}

 \def\bbZ{{\Bbb Z}}  
   \def\ZZ{\mathbb Z}

  \def\nd{\mathrel{\bigm|\kern-.7em/}}

 \def\Aut{\hbox{\rm Aut}} 
  
  \def\mod{\hbox{\rm mod}\,}

\def\Cay{\hbox{\rm Cay}} 
  
\def\qed{\hfill $\Box$} 
  \def\Arc{\hbox{\rm Arc}} 
 \def\bbZ{{\Bbb Z}}  

\def\D2n{\hbox{\rm D$_{2n}$}}

\newtheorem{thm}{Theorem}[section]
\newtheorem{cor}[thm]{Corollary}

\newtheorem{lem}[thm]{Lemma}

\makeatletter \@addtoreset{equation}{section}

\newtheorem{pro}[thm]{Proposition}

\def\pf{\noindent {\it Proof.\ }}
\def\qed{\ifmmode\square\else\nolinebreak\hfill
$\Box$\fi\par\vskip12pt}

\begin{document}

\title[Dihedral groups with the $m$-DCI property]
{Dihedral groups with the $m$-DCI property}%

\author{Jin-Hua Xie}
\address{Jin-Hua Xie, Department of Mathematics, Beijing Jiaotong University, Beijing 100044, China}
\email{jinhuaxie@bjtu.edu.cn}

\author{Yan-Quan Feng$^*$}
\address{Yan-Quan Feng, Department of Mathematics, Beijing Jiaotong University, Beijing 100044, China}
\email{yqfeng@bjtu.edu.cn}

\author{Young Soo Kwon}
\address{Young Soo Kwon, Mathematics, Yeungnam University, Kyongsan 712-749, Republic of Korea}
\email{ysookwon$@$ynu.ac.kr}

\thanks{$^*$Corresponding author}

\begin{abstract}
A Cayley digraph $\Cay(G,S)$ of a group $G$ with respect to a subset $S$ of $G$ is called {\em a CI-digraph} if for any Cayley digraph $\Cay(G,T)$ isomorphic to $\Cay(G,S)$, there is an  $\a\in \Aut(G)$ such that $S^\a=T$. For  a positive integer $m$, $G$ is said to have {\em the $m$-DCI property} if all Cayley digraphs of $G$ with out-valency $m$ are CI-digraphs. Li [The Cyclic groups with the $m$-DCI Property, European J. Combin. 18 (1997) 655-665] characterized cyclic groups with the $m$-DCI property, and in this paper, we characterize dihedral groups with the $m$-DCI property. For a dihedral group $\mathrm{D}_{2n}$  of order $2n$, assume that $\mathrm{D}_{2n}$ has the $m$-DCI property for some $1 \leq m\leq n-1$. Then it is shown that $n$ is odd, and if further $p+1\leq m\leq n-1$ for an odd prime divisor $p$ of $n$, then $p^2\nmid n$. Furthermore, if $n$ is a power of a prime $q$, then $\mathrm{D}_{2n}$ has the $m$-DCI property if and only if either $n=q$, or $q$ is odd and $1\leq m\leq q$.
\end{abstract}
\maketitle
\qquad {\textsc k}{\scriptsize \textsc {keywords.}} {\footnotesize Cayley digraph, Dihedral group, CI-digraph, DCI-group, $m$-DCI property. }

\section{Introduction}
In this paper, a {\em digraph} is an ordered pair $(V,A)$ with vertex set $V$ and arc set $A$, where $A$ is a set of ordered pairs of elements of $V$, and a digraph $(V,A)$ is called a {\em graph} if $A$ is symmetric, that is, $A=A^{-1}$, where $A^{-1}=\{(v,u)\ |\ (u,v)\in A\}$. All digraphs and graphs considered in this paper are finite and simple, and groups are finite. For a (di)graph  $\Gamma$, we use $V(\Gamma)$, $\Arc(\Gamma)$ and $\Aut(\Gamma)$ to denote the vertex set, arc set, and automorphism group of $\Gamma$, respectively. A subgroup $H$ of $\Aut(\Gamma)$  is called a {\em regular group} of automorphisms of $\Gamma$ if $H$ is regular on $V(\Gamma)$, that is, $H$ is transitive on $V(\Gamma)$ and the stabilizer $H_u=1$ for every $u\in V(\Gamma)$.

Let $G$ be a group and $S$ a subset of $G$ with $1\not\in S$. The Cayley digraph $\Cay(G,S)$ of $G$ with $S$ is defined to have vertex set $G$ and arc set $\{(g,sg)\ |\ g\in G, s\in S\}$. If $S$ is inverse-closed, that is, $S=S^{-1}$, then  $\Cay(G,S)$ is a graph.
Two Cayley (di)graphs $\Cay(G,S)$ and $\Cay(G,T)$ are called {\em Cayley isomorphic} if there is $\a \in \Aut(G)$ such that $S^\a=T$. Cayley isomorphic Cayley (di)graphs are isomorphic, but the converse is not true. A subset $S$ of $G$ with $1\not\in S$ is said to be a {\em CI-subset} if for any $\Cay(G,T)$ isomorphic to $\Cay(G,S)$,  they are Cayley isomorphic, and in  this case, $\Cay(G,S)$ is called a {\em CI-digraph}, or a {\em CI-graph} when $S=S^{-1}$. For a positive integer $m$, if all Cayley digraphs of $G$ with  out-valency $m$ are CI-digraphs, then $G$ is said to have the {\em $m$-DCI property}, and if all Cayley graphs of $G$ with valency $m$ are CI-graphs, then $G$ is said to have the {\em $m$-CI property}. Clearly, $m$-DCI property implies $m$-CI property. A group $G$ is called an {\em $m$-DCI-group} or {\em $m$-CI-group} if $G$ has the $k$-DCI property or $k$-CI property for any positive integer $k\leq m$, respectively. Furthermore, $G$ is called a {\em DCI-group} or {\em CI-group}, if $G$ have the $m$-DCI property or $m$-CI property for any positive integer $m$, respectively.

In 1967, \'Ad\'am ~\cite{Ad} conjectured that every finite cyclic group is a CI-group. Even this conjecture was disproved by Elspas and Turner~\cite{ET}, it stimulated the studies of DCI-groups and CI-groups. Classifications of cyclic DCI-groups or CI-groups were investigated by many researchers (see \cite{Alspach,Babai,Djokovic,Palfy,Turner}), and finally were completed by Muzychuk\cite{Mu1,Mu2}. There were many results on DCI-groups and CI-groups, and we refer to \cite{conder,DE1,DT,Feng,Godsil1,Kovacs,Kov,C.H.Li8,C.H.Li9,Mu3,MS,Sp,Sp1} and the references therein. Classifications of finite DCI-groups and CI-groups look very difficult, and even for dihedral groups, it is still open.  There were also some results on
$m$-DCI-groups with small $m$. For example, Fang and Xu~\cite{Fang1,Fang3,Fang2} classified abelian $m$-DCI-groups for $m\leq 3$,  and Qu and Yu~\cite{Qu} proved that for every $1\leq m\leq 3$, a dihedral group of order $2n$ ($n\geq 3$) is $m$-DCI group if and only if $n$ is odd.

A more subtle direction to study DCI-group or CI-group is to investigate the $m$-DCI property or $m$-CI property of a group. By definition, a group $G$ always has the $m$-DCI or $m$-CI property for every $m\geq |G|$, and since a (di)graph has the same automorphism group with its compliment, for each $1\leq m\leq |G|-1$, $G$ has the $m$-DCI property or $m$-CI property if and only if it has the $(|G|-m-1)$-DCI property or $(|G|-m-1)$-CI property, respectively. Therefore, we may assume that $1\leq m\leq (|G|-1)/2$. Li~\cite{C.H.Li1} characterized cyclic groups with the $m$-DCI property, and Li, Praeger and Xu~\cite{C.H.Li5} characterized all finite abelian groups with the $m$-DCI property for $1\leq m\leq 4$, where they proposed  the following problem: Characterize finite groups with the $m$-DCI property. For more details, we refer to \cite{C.H.Li0,C.H.Li2,C.H.Li3,C.H.Li4,C.H.Li6,C.H.Li7}.

In this paper, we characterize the $m$-DCI property of dihedral groups. To state the main result, we need some notations. For a finite group $G$ and a prime divisor $p$ of $|G|$, denote by $G_p$ a Sylow $p$-subgroup of $G$. For a positive integer $n$, write $\mathrm{D}_{2n}=\langle a,b\ |\ a^n=b^2=1,a^b=a^{-1}\rangle$, the dihedral group of order $2n$. Denote by $\ZZ_n$ the additive group of integers modulo $n$ and $\ZZ_n^*$ the multiplicative group of all integers coprime to $n$ in $\ZZ_n$. Clearly, $\mathrm{D}_{2}\cong\ZZ_2$ and $\mathrm{D}_{4}\cong \ZZ_2\times\ZZ_2$. By \cite{Godsil2}, $\mathrm{D}_2$ and $\mathrm{D}_4$ are DCI-groups. Note that for the $m$-DCI property of a dihedral group $\mathrm{D}_{2n}$, it suffices to consider $m$ such that $1\leq m\leq n-1$.

\begin{thm}\rm\label{mainth} Let $G$ be a dihedral group of order $2n$ with $n\geq 3$. We have the following:
\begin{itemize}
  \item [(1)] Let $G$ have the $m$-DCI property for some $1\leq m\leq n-1$. Then $n$ is odd, and if further $p+1\leq m\leq n-1$ for a prime divisor $p$ of $n$, then $G_p\cong  \bbZ_p$;
  \item [(2)] Let $n$ be a power of a prime $q$ and let $1\leq m\leq n-1$. Then $G$ has the $m$-DCI property if and only if $q$ is odd, and either $n=q$ or $1\leq m\leq q$.
\end{itemize}
\end{thm}

Based on Theorem~\ref{mainth}, we have the following corollary.

\begin{cor}\rm\label{cor}
If a dihedral group of order $2n$ ($n\geq 2$) is a $\rm{DCI}$-group, then $n=2$ or $n$ is odd-square-free.
\end{cor}

We conjecture that the converse of Corollary~\ref{cor} is true, but its proof is still elusive.

\section{Preliminaries}

In this section, we give some basic concepts and facts that will be used later. Let $\Cay(G,S)$ be a Cayley (di)graph of a group $G$ with respect to $S$. Given $g\in G$, the right multiplication $R(g): x\mapsto xg$, $\forall x\in G$, is an automorphism of $\Cay(G,S)$, and $R(G):=\{R(g)\ | \ g\in G\}$ is a regular group of automorphisms of $\Cay(G,S)$, called {\em the right regular representation} of $G$. By \cite{Babai}, we have the following well-known Babai criterion
(also see \cite[Theorem 2.4]{C.H.Li8}).

\begin{pro}\rm\label{CI-graph-prop}
A Cayley (di)graph $\Cay(G,S)$ is a CI-(di)graph if and only if every regular group of automorphisms of $\Cay(G,S)$ isomorphic to $G$ is conjugate to $R(G)$ in $\Aut(\Cay(G,S))$.
\end{pro}

The $m$-DCI property of a group is hereditary by subgroups (see~\cite[Lemma~8.2]{C.H.Li7}).

\begin{pro} \rm\label{subgroup}
Let a finite group $G$ have the $m$-DCI property for a positive integer $m$. Then every subgroup of $G$ has the $m$-DCI property.
\end{pro}

The next result gives some properties for subsets of a cyclic group (see \cite[Lemma~2.1]{C.H.Li6}).

\begin{pro} \rm\label{cyclic group property}
Let $G=\langle z\rangle$ be a cyclic group of order $n$. Assume that $i,k\in \{1,2,\ldots,n-2\}$. If $\{z,z^2,\ldots,z^k\}=\{z^i, z^{2i},\ldots.z^{ki}\}$, then $i=1$.
\end{pro}

Li~\cite[Theorem~1.2]{C.H.Li1} characterized cyclic groups with the $m$-DCI property.

\begin{pro} \rm\label{cyclicgroupMP}
Let $G$ be a cyclic group $\bbZ_n$ of order $n$ with the $m$-DCI property. For every prime divisor $p$ of $n$, if $p+1\leq m\leq n-(p+2)$ then one of the following holds:
\begin{itemize}
  \item [(1)] $G=\bbZ_{p^2}$ and $m\equiv 0$ or $-1\ (\mod p)$;
  \item [(2)] $p$ is odd and $G_p=\bbZ_p$;
  \item [(3)] $p=2$ and $G_2=\bbZ_2$ or $\bbZ_4$.
\end{itemize}
\end{pro}

A classification of dihedral DCI-groups of order $6$ times a prime was given in \cite[Theorem~1]{DE1}.

\begin{pro} \rm\label{dihedral2}  Let $p$ be a prime number. Then $\mathrm{D}_{6p}$ is a DCI-group if and only if $p\geq 5$, and in particular, $\mathrm{D}_{2p}$ is a DCI-group.
\end{pro}

From Li~\cite[Theorem~1.1]{C.H.Li2} we have the following result.

\begin{pro} \rm\label{FiniteMP}
Suppose that $G$ is an abelian group and that $p$ is the least prime divisor
of $|G|$. Then $G$ is a connected $p$-DCI-group.
\end{pro}

A finite group $G$ is called {\em homogeneous} if for any two subgroups $H$ and $K$ of $G$, every isomorphism from $H$ to $K$ can be extended to an automorphism of $G$. For dihedral groups, we have the following property (see~\cite[Lemma~1.9]{Qu}).

\begin{pro} \rm\label{homogeneous}
For an odd positive integer $n$, the dihedral group $\mathrm{D}_{2n}$ is homogeneous.
\end{pro}

\section{Sylow $p$-subgroups of $\mathrm{D}_{2n}$ with the $m$-DCI property}

Let us start with a ``general" result about isomorphic Cayley digraphs. Denote by $\overrightarrow{K}_{n,n}$ the complete bipartite digraph of order $2n$ with vertex set $\{1,2,\ldots,n\}\cup \{1',2',\ldots,n'\}$ and arc set $\{(i,j')\ |\ 1\leq i,j\leq n\}$. For a digraph $\Gamma$ and two subsets $U$ and $W$ of $V(\Gamma)$, denote by $[U,W]_\Gamma$ the  sub-digraph of
$\Gamma$ induced by $U$ and $W$, namely the vertex set is $U\cup W$ and the arc set is $\{(u,w)\in {\rm Arc}(\Gamma)\ |\ u\in U, w\in W\}$.  Let $G$ be a group and $H$ a subgroup of $G$.  A {\it right transversal} of $H$ in $G$ is a set of right coset representations of $H$ in $G$, that is, a set $L$ with $G=\cup_{\ell\in L}H\ell$ and $H\ell_1\not=H\ell_2$ for all $\ell_1,\ell_2\in L$ with $\ell_1\not=\ell_2$.

\begin{lem}\rm\label{general} Let $G$ be a finite group with $H\unlhd G$ and $H\leq M\leq G$. Assume that $C\subseteq G\backslash H$ is a union of some cosets of $H$ in $G$, and that $S_M$ and $T_M$ are subsets of $M\backslash\{1\}$ such that $S_M^\a=T_M$ for some $\a\in\Aut(M)$ and $\a$ fixes every coset of $H$ in $M$. Then $\Cay(G,S)\cong \Cay(G,T)$, where $S=C\cup S_M$ and $T=C\cup T_M$.
\end{lem}

\pf Let $\Gamma=\Cay(G,S)$ and $\Sigma=\Cay(G,T)$. Then $V(\Gamma)=V(\Sigma)=G$. Let $\{g_1,\ldots,g_\ell\}$ be a right transversal of $M$ in $G$. Define a map from $G$ to $G$ as follows:
$$\varphi: mg_k\mapsto m^\a g_k,\ \mbox{ for every } 1\leq k\leq \ell \ \mbox{ and } m\in M.$$
It is easy to see that $\varphi$ is a permutation on $G$. Clearly, $g_k^\varphi=g_k$, and since $\a$ fixes every coset of $H$ in $M$, $\varphi$ fixes every coset of $H$ in $G$. To finish the proof, it suffices to show that $\varphi$ is an isomorphism from $\Gamma$ to $\Sigma$.

Let $(v,w)$ be an arc of $\Gamma$. Then $v=mg_k\in Mg_k$ for some $1\le k\leq \ell$ and some $m\in M$. Since $S=C\cup S_M$, we have $w=gv$, where $g\in S_M$ or $C$.

Assume $g\in S_M$. Then $g\in M$. Furthermore, $v=mg_k$ and $w=gmg_k$. Thus, $v^\varphi=m^\a g_k$ and $w^\varphi=(gm)^\a g_k=g^\a m^\a g_k$ because $\a\in\Aut(M)$. Since $g^\a\in S_M^\a=T_M\subseteq T$, $(v^\varphi,w^\varphi)$ is an arc of $\Sigma$.

Assume $g\in C$. By hypothesis, $g\in Hx\subseteq C$ for some $x\in G\backslash H$. Then $g=wv^{-1} \in Hx$ and $Hg=Hx$. For any $y \in Hv$ and $z \in Hw$, $zy^{-1} \in Hwv^{-1}H = Hwv^{-1} = Hx \subseteq C$ because  $H\unlhd G$.
  So the induced sub-digraph $[Hv,Hw]_\Gamma$ is isomorphic the complete bipartite digraph $\overrightarrow{K}_{|H|,|H|}$. Since $\varphi$ fixes every coset of $H$ in $G$, we have $(Hv)^\varphi=Hv$ and $(Hw)^\varphi=Hw$, which implies $v^\varphi\in Hv$ and $w^\varphi\in Hw$. Since $g\in Hx\subseteq C\subseteq T$, similarly we have $[Hv,Hw]_\Sigma\cong \overrightarrow{K}_{|H|,|H|}$. It follows that $(v^\varphi,w^\varphi)$ is an arc of $\Sigma$.

We have proved that $(v^\varphi,w^\varphi)$ is always an arc of $\Sigma$, and hence $\varphi$ is an isomorphism from $\Gamma$ to $\Sigma$. This completes the proof. \qed

Next we give a result on isomorphic Cayley digraphs of $\mathrm{D}_{2n}$.

\begin{lem}\rm\label{n-even-con} Let $n\geq 4$ be even and let $\mathrm{D}_{2n}=\langle a,b\ |\ a^n=b^2=1,a^b=a^{-1}\rangle$. Let $H\subseteq \langle a^2\rangle$ and $K\subseteq a\langle a^2\rangle$ such that $H^{-1}=H$ and $K^{-1}=K$.  Then $\Cay(\mathrm{D}_{2n},S)\cong \Cay(\mathrm{D}_{2n},T)$, where $S=bH\cup K$ and $T=bH\cup bK$. In particular, if $\mathrm{D}_{2n}$ has the $|S|$-DCI property then $K=\emptyset$.
\end{lem}

\pf Since $n$ is even, $\langle a^2\rangle$ is a subgroup of $\langle a\rangle$ with index $2$. Then $\mathrm{D}_{2n}=\langle a^2\rangle \cup a\langle a^2\rangle \cup b\langle a^2\rangle \cup ba\langle a^2\rangle$. Define $\varphi$ by
$$\varphi: x\mapsto x \mbox{ for } x\in \langle a^2\rangle \cup b\langle a^2\rangle \mbox{ and }x\mapsto bx \mbox{ for } x\in a\langle a^2\rangle \cup ba\langle a^2\rangle,$$
that is, $\varphi$ fixes every element in $\langle a^2\rangle \cup b\langle a^2\rangle$ and interchanges $x$ and $bx$ for every $x \in a\langle a^2\rangle \cup ba\langle a^2\rangle$. Then $\varphi$ is a permutation on $\mathrm{D}_{2n}$. Set $\Gamma=\Cay(\mathrm{D}_{2n},S)$ and $\Sigma=\Cay(\mathrm{D}_{2n},T)$. We claim that $\varphi$ is an isomorphism from $\Gamma$ to $\Sigma$.

Let $(v,w)$ be an arc of $\Gamma$. Then  $w=sv$ with $s\in bH\cup K$, that is, $s\in K$ or $bH$.

\medskip
\noindent{\bf Case 1:} $s\in K\subseteq a\langle a^2\rangle$.

Since $\mathrm{D}_{2n}=\langle a^2\rangle \cup a\langle a^2\rangle \cup b\langle a^2\rangle \cup ba\langle a^2\rangle$, we have $v\in \langle a^2\rangle \cup b\langle a^2\rangle$ or $a\langle a^2\rangle \cup ba\langle a^2\rangle$. If $v\in \langle a^2\rangle \cup b\langle a^2\rangle$, then $w=sv\in a\langle a^2\rangle \cup ba\langle a^2\rangle$ as $s\in a\langle a^2\rangle$. By the definition of $\varphi$, $v^\varphi=v$ and $w^\varphi=bw=bsv$, which implies that $(v^\varphi,w^\varphi)$ is an arc of $\Sigma$ because $bs\in bK\subseteq T$. If $v\in a\langle a^2\rangle \cup ba\langle a^2\rangle$, then $w=sv\in \langle a^2\rangle \cup b\langle a^2\rangle$. It follows that $v^\varphi=bv$ and $w^\varphi=w=sv=bs^{-1}bv$ as $s\in a\langle a^2\rangle$. Since $K=K^{-1}$, we have $bs^{-1}\in bK\subseteq T$ and hence $(v^\varphi,w^\varphi)$ is an arc of $\Sigma$.

\medskip
\noindent{\bf Case 2:} $s\in bH\subseteq b\langle a^2\rangle$.

If $v\in \langle a^2\rangle \cup b\langle a^2\rangle$, then $w=sv\in \langle a^2\rangle \cup b\langle a^2\rangle$ as $s\in b\langle a^2\rangle$. Thus, $v^\varphi=v$ and $w^\varphi=w=sv$, which implies that $(v^\varphi,w^\varphi)$ is an arc of $\Sigma$ because $s\in bH\subseteq T$. If $v\in a\langle a^2\rangle \cup ba\langle a^2\rangle$, then $w=sv\in a\langle a^2\rangle \cup ba\langle a^2\rangle$. It follows that $v^\varphi=bv$ and $w^\varphi=bw=bsb(bv)$. Since $H=H^{-1}\subseteq \langle a^2\rangle$, we have $bsb\in b(bH)b=Hb=bH^{-1}=bH\subseteq T$, and hence $(v^\varphi,w^\varphi)$ is an arc of $\Sigma$.

\medskip
By Cases $1$ and $2$, $(v^\varphi,w^\varphi)$ is always an arc of $\Sigma$, and so $\varphi$ is an isomorphism from $\Gamma$ to $\Sigma$, as claimed. Thus, $\Cay(\mathrm{D}_{2n},S)\cong \Cay(\mathrm{D}_{2n},T)$.

Assume that $\mathrm{D}_{2n}$ has the $|S|$-DCI property. Since $\Cay(\mathrm{D}_{2n},S)\cong \Cay(\mathrm{D}_{2n},T)$, $\mathrm{D}_{2n}$ has an automorphism $\a$ such that $S^\a=T$. Since $n\geq 4$, $\langle a\rangle$ is characteristic in $\mathrm{D}_{2n}$, and hence $\langle a\rangle^\a=\langle a\rangle$. Suppose $K\not=\emptyset$. Since $K\subseteq S\cap \langle a\rangle$, we have $S\cap \langle a\rangle\not=\emptyset$. Thus, $T\cap \langle a\rangle=S^\a\cap \langle a\rangle^\a=(S\cap \langle a\rangle)^\a\not=\emptyset$, which is impossible as $T=bH\cup bK\subseteq b\langle a\rangle$. This completes the proof of the lemma. \qed

By \cite{Qu}, if $\mathrm{D}_{2n}$ ($n\geq 3$) has the $m$-DCI property for $m=1$ or $2$, then $n$ is odd. This is true for every $1\leq m\leq n-1$ as the following lemma.

\begin{lem}\rm\label{n-odd} Let $n\geq 3$ and let $\mathrm{D}_{2n}$ have the $m$-DCI property for some $1\leq m\leq n-1$. Then $n$ is odd.
\end{lem}

\pf Recall that $\mathrm{D}_{2n}=\langle a,b\ |\ a^n=b^2=1,a^b=a^{-1}\rangle$. Suppose to the contrary that $n$ is even. Then $a^{\frac{n}{2}}$ has order $2$. By \cite{Qu}, we may assume $m\geq 3$. If $m=3$, take $H=\{1\}$ and $K=\{a,a^{-1}\}$. By Lemma~\ref{n-even-con}, we  have $K=\emptyset$, a contradiction. Thus, we may assume $4\leq m\leq n-1$. Then $m=4k$, $4k+1$, $4k+2$ or $4k+3$, for some $k\geq 1$. Since $m\leq n-1$, we have
$2k<\frac{n}{2}$, and hence $a^i\not=a^{-i}$ for every $1\leq i\leq 2k$. Set
$$H_1=\{a^2,a^4,\ldots, a^{2k}, a^{-2},a^{-4},\ldots, a^{-2k}\},
K_1=\{a,a^3,\ldots, a^{2k-1}, a^{-1},a^{-3},\ldots, a^{-(2k-1)}\}.$$
Then $H_1=H_1^{-1}\subseteq \langle a^2\rangle$ and $K_1=K_1^{-1}\subseteq a\langle a^2\rangle$. Furthermore, $|H_1\cup K_1|=4k$.

Assume $m=4k$. Take $H=H_1$ and $K=K_1$. By Lemma~\ref{n-even-con}, $K=\emptyset$, a contradiction.

Assume $m=4k+1$. Take $H=H_1\cup\{1\}$ and $K=K_1$. Then $H=H^{-1}\subseteq \langle a^2
\rangle$. By Lemma~\ref{n-even-con}, $K=\emptyset$, a contradiction.

Assume $m=4k+2$. If $\frac{n}{2}$ is odd then take $H=H_1\cup\{1\}$ and $K=K_1\cup\{a^{\frac{n}{2}}\}$, and if $\frac{n}{2}$ is even, take $H=H_1\cup\{1,a^{\frac{n}{2}}\}$ and $K=K_1$. Then $H=H^{-1}\subseteq \langle a^2
\rangle$ and $K=K^{-1}\subseteq a\langle a^2
\rangle$. By Lemma~\ref{n-even-con}, $K=\emptyset$, a contradiction.

Assume $m=4k+3$. Since $1\leq m\leq n-1$, we have $1\leq 2k\leq \frac{n-4}{2}=\frac{n}{2}-2<\frac{n}{2}-1$.
If $\frac{n}{2}-1$ is even then take $H=H_1\cup\{1,a^{\frac{n}{2}-1},a^{-(\frac{n}{2}-1)}\}$ and $K=K_1$, and if $\frac{n}{2}-1$ is odd, take $H=H_1\cup\{1\}$ and $K=K_1\cup \{a^{\frac{n}{2}-1},a^{-(\frac{n}{2}-1)}\}$. Then $H=H^{-1}\subseteq \langle a^2
\rangle$ and $K=K^{-1}\subseteq a\langle a^2
\rangle$. By Lemma~\ref{n-even-con}, $K=\emptyset$, a contradiction. \qed

Now we consider Sylow $p$-subgroups of the dihedral group $\mathrm{D}_{2n}$ with $m$-DCI property, where $p+1\leq m\leq n-1$ for some prime divisor $p$ of $2n$.

\begin{lem}\rm\label{p-odd}  Let $G=\mathrm{D}_{2n}$ ($n\geq 3$) have the $m$-DCI property. If $p$ is a prime divisor of $2n$ such that $p+1\leq m\leq n-1$, then $G_p\cong  \ZZ_p$.
\end{lem}

\pf Let $p$ be a prime divisor of $2n$ such that $p+1\leq m\leq n-1$. By Lemma~\ref{n-odd}, $G_2\cong  \ZZ_2$. Thus, we may assume that $p$ is odd. Note that $G=\mathrm{D}_{2n}=\langle a,b\ |\ a^n=b^2=1,a^b=a^{-1}\rangle$.
Let
$P=\langle z\rangle$ with $z=a^{n/p}$. Then $P$ is the unique subgroup of order $p$ in $G$, and hence characteristic in $G$. For convenience, write $n'=n/p$ and $n_p=p^d$, where $n_p$ is the largest $p$-power dividing $n$. In particular, $z=a^{n'}$. Similarly, denote by $n_{p'}$ the largest divisor of $n$ such that $p\nmid n_{p'}$. Then $(n_p,n_{p'})=1$ and $n=n_pn_{p'}$.

Suppose, for a contradiction,  that $G_p\not\cong \ZZ_p$.  Then $d\geq 2$ and $p\mid n'$. By Proposition~\ref{subgroup}, $\langle a\rangle$ has the $m$-DCI property. Let us first consider the case for $m\not\equiv 0, -1\ (\mod p)$.

\medskip
\noindent{\bf Case 1:} $m\not\equiv 0, -1\ (\mod p)$

Since $\langle a\rangle$ has the $m$-DCI property, if $p+1\leq m\leq n-(p+2)$ then Proposition~\ref{cyclicgroupMP} implies $\langle a\rangle_p\cong \ZZ_p$, where $\langle a\rangle_p$ is the Sylow $p$-subgroup of $\langle a\rangle$. It follows that $G_p=\langle a\rangle_p\cong\ZZ_p$, a contradiction. Thus, $n-(p+1)\leq m\leq n-1$ as $p+1\leq m\leq n-1$. Note that $d\geq 2$ implies $p^2\mid n$.
Since $m\not \equiv 0, -1\ (\mod p)$, we have  $n-p+1\leq m\leq n-2=n-p+p-2$, and so
$m=kp+j$, where $k=(n-p)/p=n'-1$ and $1\leq j\leq p-2$. Let
\[\begin{cases}
S=\{a,ba,\ldots,ba^{k-1}\}P \cup \{z,\ldots,z^j\} ;\\
T=\{a,ba,\ldots,ba^{k-1}\}P \cup \{z^{-1},\ldots,z^{-j}\}.
\end{cases} \]
Here $\{a,ba,\ldots,ba^{k-1}\}P:=aP\cup baP\cup\ldots\cup ba^{k-1}P$. Since $P$ is abelian, it has an automorphism $\alpha$ such that $z^\alpha=z^{-1}$. Thus, $\{z,\ldots,z^j\}^\alpha=\{z^{-1},\ldots,z^{-j}\}$. Taking $M=H=P$ and $C=\{a,ba,\ldots,ba^{k-1}\}P$, by Lemma~\ref{general} we have $\Cay(G,S)\cong \Cay(G,T)$. Since $G$ has the $m$-DCI property, there is $\sigma\in \Aut(G)$ such that with $S^\sigma=T$. Then $P^\sigma=P$ as $P$ is characteristic in $G$.

Since $\langle a\rangle$ is characteristic in $G$, we have $a^\sigma=a^r$ for some $r\in\ZZ_n^*$. Let $x\in aP$. Then $x=az^l=a^{ln'+1}$ for some $0\leq l\leq p-1$. Then $x$ has order $n/(n,ln'+1)$. Since $p\mid n'$, we have $(p,ln'+1)=1$. Since  $(n')_{p'}=n_{p'}$, we have $(n_{p'}, ln'+1)=1$, and hence  $(n,ln'+1)=1$. It follows $o(x)=n$, that is, every element in $aP$ has order $n$. Furthermore, every element in $\{ba,\ldots,ba^{k-1}\}P$ has order $2$, and every element in $\{z,\ldots,z^j\}$ has order $p$. Since $\sigma\in \Aut(G)$, we have $(aP)^\sigma=aP$ and $\{z,\ldots,z^j\}^\sigma=\{z^{-1},\ldots,z^{-j}\}$. From the latter, $\{z,\ldots,z^j\}=\{z^{-r},\ldots,z^{-rj}\}$, and by Proposition~\ref{cyclic group property},
$r=-1\ (\mod p)$, which implies $r=sp-1$ for some integer $s$. From the former, $a^{sp-1}P=aP$ and hence $a^{sp-2}\in P$, which is impossible as $p^d\mid o(a^{sp-2})$.

\medskip
\noindent{\bf Case 2:} $m\equiv -1\, (\mod p)$

We divide this case into two subcases: $n\not=p^2$ and $n=p^2$.

\medskip
\noindent{\bf Subcase 2.1:} $n\not=p^2$.

If $p+1\leq m\leq n-(p+2)$, Proposition~\ref{cyclicgroupMP} implies $\langle a\rangle_p\cong \ZZ_p$ and hence $G_p\cong \ZZ_p$, contradicting that $d\geq 2$. Thus, we may assume $n-(p+1)\leq m\leq n-1$, and since $m\equiv -1\ (\mod p)$, we have $m=n-2p+(p-1)$ or $n-p+(p-1)$. It follows that $m=kp+(p-1)$ with $k=n'-2$ or $n'-1$.

First assume that $n$ is not a $p$-power. Then $n=n_pn_{p'}$ with $n_{p'}\not=1$ and let
\[\begin{cases}
S=\{a,ba,\ldots,ba^{k-1}\}P\cup \{z,\ldots,z^{p-2}, a^{n_p}\},\\
T=\{a,ba,\ldots,ba^{k-1}\}P\cup \{z^{-1},\ldots,z^{-(p-2)}, a^{n_p}\};
\end{cases} \]

Let $M=\langle z,\ldots,z^{p-2}, a^{n_p}\rangle$. Note that $o(a^{n_p})=n/n_p=n_{p'}$ and $M$ is a cyclic group of order $pn_{p'}$ of $\langle a\rangle$. Thus, $M$ is characteristic in $G$, and since $(p,n_{p'})=1$, $M$ has an automorphism $\alpha$ such that $z^\alpha=z^{-1}$ and $(a^{n_p})^\alpha=a^{n_p}$. It follows that $\{z,\ldots,z^{p-2}, a^{n_p}\}^\a=\{z^{-1},\ldots,z^{-(p-2)}, a^{n_p}\}$ and $\a$ fixes every coset of $P$ in
$M$. By Lemma~\ref{general}, $\Cay(G,S)\cong \Cay(G,T)$, and so there exists $\sigma\in \Aut(G)$ such that $S^\sigma=T$. In particular, $a^\sigma=a^r$ for some $r\in\ZZ_n^*$ and $P^\sigma=P$.

Considering the orders of elements in $S$ and $T$, we have $(aP)^\sigma=aP$ and $\{z,\ldots,z^{p-2}\}^\sigma=\{z^{-1},\ldots,z^{-(p-2)}\}$, that is, $\{z^{-r},\ldots,z^{-r(p-2)}\}=\{z,\ldots,z^{p-2}\}$, which implies $r=-1\ (\mod p)$ by Proposition~\ref{cyclic group property}. Then $r=sp-1$ and $(aP)^\sigma=aP$ implies $a^{sp-2}P=P$, which is impossible because $p^d\mid o(a^{sp-2})$.

\medskip

Now assume that $n$ is a power of $p$. Since $n\not=p^2$, we have $n=p^d$ with $d\geq 3$. Furthermore, $n'=p^{d-1}$, $m=kp+(p-1)$ with $k=p^{d-1}-2$ or $p^{d-1}-1$. Let
\[\begin{cases}
S=\{a,ba,\ldots,ba^{k-1}\}P\cup \{z,\ldots,z^{p-2}, a^{p^{d-2}}\},\\
T=\{a,ba,\ldots,ba^{k-1}\}P\cup \{z,\ldots,z^{p-2}, a^{p^{d-2}+p^{d-1}}\};
\end{cases} \]

Note that $P=\langle a^{p^{d-1}}\rangle$ and $ a^{p^{d-2}}P= a^{p^{d-2}+p^{d-1}}P$. Every element in $a^{p^{d-2}}P$ has order $p^2$. Write $M=\langle z,\ldots,z^{p-2}, a^{p^{d-2}}\rangle$. Then $M=\langle z,\ldots,z^{p-2},  a^{p^{d-2}+p^{d-1}}\rangle=\langle a^{p^{d-2}} \rangle=\langle a^{p^{d-2}+p^{d-1}}\rangle$, which has order $p^2$. Thus, $M$ has an automorphism $\a$ induced by $a^{p^{d-2}}\mapsto a^{p^{d-2}+p^{d-1}}$. It follows that $\{z,\ldots,z^{p-2}, a^{p^{d-2}}\}^\a=\{z,\ldots,z^{p-2}, a^{p^{d-2}+p^{d-1}}\}$ and $\a$ fixes every coset of $P$ in $M$. By Lemma~\ref{general}, $\Cay(G,S)\cong \Cay(G,T)$, and so there exists $\sigma\in \Aut(G)$ such that $S^\sigma=T$. In particular, $a^\sigma=a^r$ for some $r\in\ZZ_{p^d}^*$.

Considering the orders of elements in $S$ and $T$, we have $(aP)^\sigma=aP$ and $(a^{p^{d-2}})^\sigma= a^{p^{d-2}+p^{d-1}}$. Then $\sigma\not=1$. Let $a^\sigma=a^r$ for some $r\in\ZZ_{p^d}^*$. Then $a^{rp^{d-2}}=(a^{p^{d-2}})^\sigma= a^{p^{d-2}+p^{d-1}}$, and $p^d\mid ((r-1)p^{d-2}-p^{d-1})$. It follows that $r=\ell p+1$ for some integer $\ell$, and hence  $a^{\ell p^{d-1}+p^{d-2}} =a^{rp^{d-2}}= a^{p^{d-2}+p^{d-1}}$, that is, $a^{\ell p^{d-1}} =a^{p^{d-1}}$. This implies that $(\ell,p)=1$. Since $(aP)^\sigma=aP$, we have $a^{\ell p}P=P$ and so  $a^{\ell p}\in P$, which is impossible because $(\ell,p)=1$ implies that $a^{\ell p}$ has order $p^{d-1}\geq p^2$.

\medskip
\noindent{\bf Subcase 2.2:} $n=p^2$.

In this case, $n'=p$ and $z=a^p$. Since $p+1\leq m\leq n-1=p^2-1$ and $m\equiv -1\ (\mod p)$, we have $m=kp+(p-1)$, where $1\leq k\leq p-1$. Let
\[\begin{cases}
S=\{a,ba,\ldots,ba^{k-1}\}P \cup \{z,\ldots,z^{p-2}, b\};\\
T=\{a,ba,\ldots,ba^{k-1}\}P \cup \{z^{-1},\ldots,z^{-(p-2)}, b\}.
\end{cases} \]

Let $M=\langle z,\ldots,z^{p-2}, b\rangle=\langle a^p,b\rangle$. Then $M$ is a dihedral group of order $2p$, and has an automorphism $\a$ induced by $a^p\mapsto a^{-p}$ and $b\mapsto b$. Furthermore, $\{z,\ldots,z^{p-2}, b\}^\a=\{z^{-1},\ldots,z^{-(p-2)}, b\}$ and $\a$ fixes every coset of $P$ in $M$. By Lemma~\ref{general}, $\Cay(G,S)\cong \Cay(G,T)$, and so there exists $\sigma\in \Aut(G)$ such that $S^\sigma=T$. Then $a^\sigma=a^r$ for some $r\in\ZZ_{p^2}^*$ and $P^\sigma=P$.

Considering the orders of elements in $S$ and $T$, we have $(aP)^\sigma=aP$ and $\{z,\ldots,z^{p-2}\}^\sigma=\{z^{-1},\ldots,z^{-(p-2)}\}$, which implies that $r=sp-1$.
Then  $(aP)^\sigma=aP$ implies $a^{sp-2}P=P$, which is impossible because $o(a^{sp-2})=p^2$.

\medskip
\noindent{\bf Case 3:} $m\equiv 0\, (\mod p)$.

We also divide this case into two subcases: $n\not=p^2$ and $n=p^2$.

\medskip
\noindent{\bf Subcase 3.1:} $n\not=p^2$.

If $p+1\leq m\leq n-(p+2)$, Proposition~\ref{cyclicgroupMP} implies $\langle a\rangle_p\cong \ZZ_p$  and then $G_p\cong \ZZ_p$, a contradiction. Thus we may assume $n-(p+1)\leq m\leq n-1$. Since $m\equiv 0\ (\mod p)$, we have $m=n-p=kp$ with $k=n'-1$.

Assume that $n$ is not a $p$-power. Then $n=n_pn_{p'}$ with $n_{p'}\not=1$ and let
\[\begin{cases}
S=\{a,ba,\ldots,ba^{k-2}\}P\cup \{z,\ldots,z^{p-2}, a^{n_p}, b\},\\
T=\{a,ba,\ldots,ba^{k-2}\}P\cup \{z^{-1},\ldots,z^{-(p-2)}, a^{n_p}, b\};
\end{cases} \]

Let $M=\langle z,\ldots,z^{p-2}, a^{n_p},b\rangle$. Since $o(a^{n_p})=n/n_p=n_{p'}$, $M$ is a dihedral group of order $2pn_{p'}$, and has an automorphism $\alpha$ such that $z^\alpha=z^{-1}$, $(a^{n_p})^\alpha=a^{n_p}$ and $b^\a=b$. It follows that $\{z,\ldots,z^{p-2}, a^{n_p}, b\}^\a=\{z^{-1},\ldots,z^{-(p-2)}, a^{n_p}, b\}$ and $\a$ fixes every coset of $P$ in $M$.
By Lemma~\ref{general}, $\Cay(G,S)\cong \Cay(G,T)$, and so there exists $\sigma\in \Aut(G)$ such that $S^\sigma=T$. In particular, $a^\sigma=a^r$ for some $r\in\ZZ_n^*$ and $P^\sigma=P$. Considering the orders of elements in $S$ and $T$, we have $(aP)^\sigma=aP$ and $\{z,\ldots,z^{p-2}\}^\sigma=\{z^{-1},\ldots,z^{-(p-2)}\}$, that is, $\{z^{-r},\ldots,z^{-r(p-2)}\}=\{z,\ldots,z^{p-2}\}$, which implies $r=-1\ (\mod p)$ by Proposition~\ref{cyclic group property}. Thus $r=sp-1$, and then $(aP)^\sigma=aP$ implies $a^{sp-2}P=P$, which is impossible because $p^d\mid o(a^{sp-2})$.

\medskip

Now assume that $n$ is a power of $p$. Then $n=p^d$ with $d\geq 3$ and $n'=p^{d-1}$. Let
\[\begin{cases}
S=\{a,ba,\ldots,ba^{k-2}\}P\cup \{z,\ldots,z^{p-2}, a^{p^{d-2}},b\},\\
T=\{a,ba,\ldots,ba^{k-2}\}P\cup \{z,\ldots,z^{p-2}, a^{p^{d-2}+p^{d-1}},b\};
\end{cases} \]

Since $P=\langle a^{p^{d-1}}\rangle$, we have $a^{p^{d-2}}P= a^{p^{d-2}+p^{d-1}}P$. Let $M=\langle z,\ldots,z^{p-2}, a^{p^{d-2}},b\rangle$. Since $\langle z,\ldots,z^{p-2},  a^{p^{d-2}+p^{d-1}}\rangle=\langle a^{p^{d-2}} \rangle=\langle a^{p^{d-2}+p^{d-1}}\rangle$, $M$ is a dihedral group of order $2p^2$. Then $M$ has an automorphism $\a$ induced by $a^{p^{d-2}}\mapsto a^{p^{d-2}+p^{d-1}}$ and $b\mapsto b$, which fixes every coset of $P$ in $M$. Furthermore, $\{z,\ldots,z^{p-2}, a^{p^{d-2}},b\}^\a=\{z,\ldots,z^{p-2}, a^{p^{d-2}+p^{d-1}},b\}$. By Lemma~\ref{general}, $\Cay(G,S)\cong \Cay(G,T)$, and so there exists $\sigma\in \Aut(G)$ such that $S^\sigma=T$. In particular, $a^\sigma=a^r$ for some $r\in\ZZ_{p^d}^*$.

Considering the orders of elements in $S$ and $T$, we have $(aP)^\sigma=aP$ and $(a^{p^{d-2}})^\sigma= a^{p^{d-2}+p^{d-1}}$. Then $\sigma\not=1$ and $a^{rp^{d-2}}=(a^{p^{d-2}})^\sigma= a^{p^{d-2}+p^{d-1}}$, forcing $p^d\mid ((r-1)p^{d-2}-p^{d-1})$. It follows that $r=\ell p+1$ for some integer $\ell$, and hence  $a^{\ell p^{d-1}+p^{d-2}} =a^{rp^{d-2}}= a^{p^{d-2}+p^{d-1}}$, that is, $a^{\ell p^{d-1}} =a^{p^{d-1}}$. This implies that $(\ell,p)=1$. Since $(aP)^\sigma=aP$, we have $a^{\ell p}P=P$ and so  $a^{\ell p}\in P$, which is impossible because $(\ell,p)=1$ implies that $a^{\ell p}$ has order $p^{d-1}\geq p^2$.

\medskip
\noindent{\bf Subcase 3.2:} $n=p^2$.

In this case, $n'=p$ and $z=a^p$. Since $p+1\leq m\leq n-1=p^2-1$ and $m\equiv 0\ (\mod p)$, we have $m=kp$, where $2\leq k\leq p-1$. Let
\[\begin{cases}
S=\{a,ba,\ldots,ba^{k-2}\}P \cup \{z,\ldots,z^{p-2}, b,bz\};\\
T=\{a,ba,\ldots,ba^{k-2}\}P \cup \{z^{-1},\ldots,z^{-(p-2)}, b,bz^{-1}\}.
\end{cases} \]
Let $M=\langle z,b\rangle=\langle a^p,b\rangle$. Then $M$ is a dihedral group of order $2p$, and has an automorphism $\a$ induced by $a^p\mapsto a^{-p}$ and $b\mapsto b$. It is easy to see that $\{z,\ldots,z^{p-2}, b,bz\}^\a=\{z^{-1},\ldots,z^{-(p-2)}, b,bz^{-1}\}$ and $\a$ fixes every coset of $P$ in $M$. By Lemma~\ref{general}, $\Cay(G,S)\cong \Cay(G,T)$. Then there exists $\sigma\in \Aut(G)$ such that $S^\sigma=T$. In particular, $a^\sigma=a^r$ for some $r\in\ZZ_{p^2}^*$.

Note that every element in $aP$, $\{ba,\ldots,ba^{k-2}\}P\cup\{b,bz,bz^{-1}\}$, or $\{z,\ldots,z^{p-2}, z^{-1},\ldots$, $z^{-(p-2)}\}$ has order $p^2$, $2$ or $p$, respectively. We have $(aP)^\sigma=aP$ and $\{z,\ldots,z^{p-2}\}^\sigma=\{z^{-1},\ldots,z^{-(p-2)}\}$. The latter implies $r=sp-1$, and then the former implies $a^{sp-2}P=P$, which is impossible because $o(a^{sp-2})=p^2$. \qed

\section{Proof of Theorem~\ref{mainth}\label{proof}}
In this section, we prove Theorem~\ref{mainth}. The following lemma shows that a dihedral group has ``similar" properties on CI-subset with its unique cyclic subgroup.

\begin{lem}\rm\label{conjugate} Let $\mathrm{D}_{2n}=\langle a,b\ |\ a^n=b^2=1, a^b=a^{-1}\rangle$ with $n\geq 3$, and let $\Gamma=\Cay(\mathrm{D}_{2n},S)$ be a Cayley digraph of $\mathrm{D}_{2n}$ with $A=\Aut(\Gamma)$. Assume that $A$ has a regular subgroup $R$ isomorphic $\mathrm{D}_{2n}$. Then $R$ and $R(\mathrm{D}_{2n})$ are conjugate in $A$ if and only if the unique cyclic subgroup of order $n$ in $R$ is conjugate to $\langle R(a)\rangle$ in $A$.
\end{lem}

\pf The necessity is clear because a dihedral group of order $2n$ ($n\geq 3$) has a unique cyclic subgroup of order $n$. To prove the sufficiency, we may write $R=\langle x\rangle\rtimes\langle y\rangle$ with $o(x)=n$, $o(y)=2$ and $x^y=x^{-1}$, and assume that $A$ has an element, say $\a$, such that $\langle x\rangle^\a=\langle R(a)\rangle$. To finish the proof, it suffices to show that
$R^\a=R(\mathrm{D}_{2n})$.

Clearly, $R^\a=\langle x\rangle^\a\rtimes\langle y\rangle^\a =\langle R(a)\rangle\rtimes\langle y^\a\rangle \cong \mathrm{D}_{2n}$ is regular on $V(\Gamma)$ and $R(a)^{y^\a}=R(a)^{-1}$. Since $R(a)^{R(b)}=R(a)^{-1}$, we have $R(a)^{y^\a R(b)}=R(a)$. Note that $\langle R(a)\rangle$ has exactly two orbits on $V(\Gamma)$, that is, $\langle a\rangle$ and $b\langle a\rangle$. Since both $R^\a$ and $R(\mathrm{D}_{2n})$ are regular on $V(\Gamma)$ and contain the normal subgroup $\langle R(a)\rangle$, both $R(b)$ and $y^\a$ interchanges $\langle a\rangle$ and $b\langle a\rangle$, respectively. In particular, $R(b)$ interchanges $1$ and $b$, and $y^\a$ interchanges $1$ and $ba^i$ for some $i\in\ZZ_n$. Let $z=y^\a R(a)^{-i}\in R^\a$. Then $1^z=b$, and $z$ interchanges $\langle a\rangle$ and $b\langle a\rangle$. It follows that  $R(b)z$ fixes  $\langle a\rangle$ and $b\langle a\rangle$ setwise, respectively.

Since $\langle R(a)\rangle \rtimes\langle {y^\a}\rangle\cong \mathrm{D}_{2n}$, $z=y^\a R(a)^{-i}$ is an involution and hence $b^z=1$ as $1^z=b$. Thus, $R(b)z$ fixes $1$ and $b$, and so $\langle R(a),R(b)z\rangle$ has two orbits on $V(\Ga)$, that is, $\langle a\rangle$ and $b\langle a\rangle$. Since $R(a)^z=R(a)^{y^\a R(a)^{-i}}=R(a)^{-1}$, we have $R(a)^{R(b)z}=R(a)$ and so $\langle R(a),R(b)z\rangle$ is abelian. This implies that $R(b)z$ fixes $\langle a\rangle$ and $b\langle a\rangle$ pointwise. It follows that $R(b)z=1$, that is $R(b)=z$, and so $R^\a= \langle R(a)\rangle\rtimes\langle y^\a\rangle =\langle R(a)\rangle\rtimes\langle y^\a R(a)^{-i}\rangle=\langle R(a)\rangle\rtimes\langle z\rangle=\langle R(a)\rangle\rtimes\langle R(b)\rangle=R(\mathrm{D}_{2n})$. This completes the proof. \qed


Next we give a  sufficiency condition for CI-digraphs on dihedral groups.

\begin{lem}\rm\label{p-CI} Let $\Gamma=\Cay(\mathrm{D}_{2n},S)$ be a digraph with $n$ odd and let $A=\Aut(\Gamma)$. If $(|A_1|,n)=1$ then $\Gamma$ is a CI-digraph, and if $\Gamma$ is connected and $|S|<p$ for the least prime divisor $p$ of $n$, then $\Gamma$ is a CI-digraph.
\end{lem}

\pf Let $\mathrm{D}_{2n}=\langle a,b\ |\ a^n=b^2=1, a^b=a^{-1}\rangle$ and let $\pi=\{q\ |\  q \mbox{ is  a prime divisor of } n\}$. Assume that $(|A_1|,n)=1$. Since $A=R(G)A_1$ and $n$ is odd, $\langle R(a)\rangle$ is a Hall $\pi$-subgroup of $A$. By \cite[Theorem 9.1.10]{Robinson}, all nilpotent Hall $\pi$-subgroup of $A$ are conjugate, and then Lemma~\ref{conjugate} implies that all regular subgroups of $A$ isomorphic to $R(G)$ are conjugate. By Proposition~\ref{CI-graph-prop}, $\Gamma$ is a CI-digraph.

Now assume that $\Gamma$ is connected and $|S|<p$, where $p$ is the least prime divisor of $n$. Then $\Gamma$ has out-valency less than $p$, and the connectivity of $\Gamma$ implies that each prime divisor of $|A_1|$ is less than $p$. It follows that $(|A_1|,n)=1$ and then the above paragraph means that $\Gamma$ is a CI-digraph. \qed

Now we are ready to prove Theorem~\ref{mainth}.
\medskip

\noindent{\bf Proof of Theorem~\ref{mainth}:} Let $G=\mathrm{D}_{2n}=\langle a,b\ |\ a^n=b^2=1,a^b=a^{-1}\rangle$ be the dihedral group of order $2n$ with $n\geq 3$.

To prove part~(1), assume that $G$ has the $m$-DCI property for some $1\leq m\leq n-1$, and that $p$ is a prime divisor of $2n$. Then $n$ is odd by Lemma~\ref{n-odd}. If further $p+1\leq m\leq n-1$. then  $G_p\cong \bbZ_p$  by  Lemma~\ref{p-odd}. This completes the proof of part~(1).

To prove part~(2), assume $n=q^r$ with $q$ a prime and $r$ a positive integer, and assume  $1\leq m\leq n-1$.

To prove the necessity of part~(2), let $G$ be of the $m$-DCI property. By part~(1), $n$ is odd, and hence $q$ is odd. Now assume that $r\geq 2$. If $q+1\leq m\leq n-1$, then $G_q\cong\ZZ_q$ by part~(1), contradicting $r\geq 2$. Thus $1\leq m \leq q$, completing the proof of the necessity.

To prove the sufficiency of part~(2), assume that $q$ is odd, and that either $n=q$ or $1\leq m\leq q$. For $n=q$, $\mathrm{D}_{2q}$ is a DCI-group by Proposition~\ref{dihedral2}, and hence has the $m$-DCI property. Thus, we may assume $1\leq m\leq q$. Let $\Gamma=\Cay(G,S)$ with $|S|=m$ and let $\Sigma=\Cay(G,T)$ such that $\Gamma\cong\Sigma$. To finish the proof, we only need to show that $G$ has an automorphism $\a$ such that $S^\a=T$.

Let $K=\langle S\rangle$ and $H=\langle T\rangle$. Then $K\leq G$ and $H\leq G$. Note that $\Gamma=|G:K|\,\Cay(K,S)$, the  $|G:K|$ distinct copies of the connected digraph $\Cay(K,S)$, and similarly, $\Sigma=|G:H|\,\Cay(H,T)$, the $|G:H|$ distinct copies of the connected digraph $\Cay(H,T)$. Since $\Sigma\cong\Gamma$, we have $|K|=|H|$, and since $G$ is the dihedral group of order $2q^r$ with $q$ odd, both $K$ and $H$ are either a cyclic subgroup of $\langle a\rangle$ or a dihedral subgroup of $G$, with the same order, which implies $K\cong H$. By Proposition~\ref{homogeneous}, there is $\b\in\Aut(G)$ such that $K^\b=H$, and then $\Cay(K,S)\cong \Cay(H,S^\b)$. It follows that $\Cay(H,S^\b)\cong \Cay(H,T)$. Let $\Delta=\Cay(H,S^\b)$ and $B=\Aut(\Delta)$.

Now we claim that $\Delta$ is a CI-digraph. Recall that $H$ is cyclic subgroup of $\langle a\rangle$ or a dihedral subgroup of $G$. If $H\leq \langle a\rangle$ then Proposition~\ref{FiniteMP} implies that $H$ is a $q$-DCI-group and hence has the $m$-DCI property for each $1\leq m\leq q$. So  $\Delta$ is a CI-digraph. Now we may assume that  $H$ is dihedral. If $q\nmid |B_1|$ (in this case $(|B_1|,q^r)=1$) or $|S^\b|<q$, Lemma~\ref{p-CI} implies that $\Delta$ is a CI-digraph. Thus, we may further assume that $q\mid |B_1|$ and $|S^\b|=q$, which implies that $\Delta$ is arc-transitive. Since $H$ is dihedral and $\Delta$ is connected, $S^\b$ contains an  involution, and then the arc-transitivity of $\Delta$ implies that $\Delta$ is a graph, that is, $S^\b=(S^\b)^{-1}$. The classification of connected arc-transitive Cayley graphs of a dihedral group of order $2$ times $q$-power, given by Kov\'{a}cs~\cite{Kovacs}, implies that  $\Delta\cong K_{q,q}$. It follows that $H\cong \mathrm{D}_{2q}$, and since $\mathrm{D}_{2q}$ is a DCI-group, $\Delta$ is a CI-digraph as claimed.

Since $\Cay(H,S^\b)\cong \Cay(H,T)$ and $\Delta$ is a CI-digraph, $H$ has an automorphism $\g'$ such that $(S^\b)^{\g'}=T$. By Proposition~\ref{homogeneous}, $G$ has an automorphism $\g$ such that $\g|_H=\g'$, that is, the restriction of $\g$ on $H$ is $\g'$. It follows that $S^{\b\g}=T$. Take $\a=\b\g$. Then $\a\in\Aut(G)$ and $S^\a=T$. This completes the proof.  \qed

\medskip
\noindent {\bf Acknowledgements:} The work was partially supported by the National Natural Science Foundation of China (12161141005,12271024) and the 111 Project of China (B16002).  The third author was partially supported by the Basic Science Research Program through the National Research Foundation of Korea (NRF)
funded by the Ministry of Education (2018R1D1A1B05048450) and (2021K2A9A2A11101586).


\end{document}